\documentclass{article}
\usepackage{amsmath,amsfonts, amssymb}

\newcommand{\E}{\mathbf{E}}
\newcommand{\R}{\mathbf{R}}
\newcommand{\p}{\mathbf{P}}

\newtheorem{lemma}{Lemma}
\newtheorem{thm}{Theorem}

\newtheorem{prop}[thm]{Proposition}
\newtheorem{defn}{Definition}

  \newenvironment{Proof}{\noindent{\bf Proof} \ }{\QED}\smallskip
\newcommand\QED{\newline \rightline{$\blacksquare$} \bigskip}

\title{Classifying the Concentration of the Boolean Cube for Dependent Distributions}
\author{Jonathan Root\footnote{ Department of Mathematics and Statistics, Boston University, Boston, MA 02215, USA}, 
Mark Kon \footnote{ Department of Mathematics and Statistics, Boston University, Boston, MA 02215, USA}
}
\begin{document}

\maketitle
\begin{abstract}
A metric probability space $(\Omega,d)$ obeys the {\it concentration of measure phenomenon} if subsets of measure $1/2$ enlarge to subsets of measure close to 1 as a transition parameter $\epsilon$ approaches a limit. 
In this paper we consider the concentration of the space itself, namely the concentration of the metric $d(x,y)$ for a fixed $y\in \Omega$.  For any $y\in \Omega$, the concentration of $d(x,y)$ is guaranteed for product distributions in high dimensions $n$, as $d(x,y)$ is a Lipschitz function in $x$. In fact, in the product setting, the rate at which the metric concentrates is of the same order in $n$ for any fixed $y\in \Omega$. The same thing, however, cannot be said for certain dependent (non-product) distributions. For the Boolean cube $I_n$ (a widely analyzed simple model), we show that, for any dependent distribution, the rate of concentration of the Hamming distance $d_H(x,y)$, for a fixed $y$, depends on the choice of $y\in I_n$, and on the variance of the conditional distributions $\mu(x_k \mid x_1,\dots, x_{k-1})$, $2\leq k\leq n$. We give an inductive bound which holds for all probability distributions on the Boolean cube, and characterize the quality of concentration by a certain positive (negative) correlation condition.
Our method of proof is advantageous in that it is both simple and comprehensive. We consider uniform bounding techniques when the variance of the conditional distributions is negligible, and show how this basic technique applies to the concentration of the entire class of Lipschitz functions on the Boolean cube. 
\end{abstract}

\section{Introduction}

The concentration of measure phenomenon as it is now known
originated with Milman in his proof of Dvoretzky's theorem \cite{milman,asymp, barvinok}.
 The so-called Gromov-Milman formulation \cite{gromov-milman,asymp} of the concentration
of measure phenomenon begins with a (Polish) metric space
$(\Omega,d)$ provided with a Borel probability measure $\mu$. 
The ideas are based on the principle that to 
develop good tail bounds on the deviation of a general
Lipschitz function on $\Omega$ about its median, one needs only
prove that sets of measure at least 1/2 enlarge to (are metrically close to) sets of measure close
to one. Specifically,
the ($\epsilon-$) enlargement of a set $A$ is defined by $A_{\epsilon} = \{x\in \Omega : d(x,A)\leq \epsilon\}$.
One may then define the {\it concentration function} \cite{asymp}
$\alpha(\Omega,\epsilon)$, for any $\epsilon>0$, by 
\begin{equation}\label{concentrationfunction}
\alpha(\Omega,\epsilon) = 1-\inf\{\mu(A_{\epsilon}) : A\subset \Omega \text{ Borel with } \mu(A)\geq 1/2\}.
\end{equation}
With respect to a family $(\Omega_n, d_n,\mu_n)$, $n\geq 1$ of metric probability
spaces, we have the notion of concentration about one point if for every $\epsilon>0$,
\begin{equation}\label{concentrationdef}
\alpha(\Omega_n, \epsilon\text{ diam }\Omega_n) \to 0, \; (n\to \infty),
\end{equation}
and we have {\it normal concentration} if 
\begin{equation}\label{normalconcentration}
\alpha(\Omega_n, \epsilon) \leq c_1\exp(-c_2\epsilon^2 n),
\end{equation}
for some constants $c_1,c_2$. Typically bounds like (\ref{normalconcentration}) describe the tail 
behavior (in the variable $x\in \Omega_n=\R^n$) of an $n$-dimensional probability distribution as $n$ becomes large. 
This is perhaps the most general formulation of the concentration of measure phenomenon, but a specialization of this definition is of great importance in modern probability theory and its applications. 

Suppose we are
given a continuous real-valued function $f\in C(\Omega_n)$. Denote its {\it modulus of continuity} by
$\omega_f(\epsilon) = \sup\{|f(x)-f(y)| : d(x,y) \leq \epsilon\}$. We will assume the function is ``well-behaved," meaning that $\omega_f(\epsilon)$ is small.
Let $M_f$ denote the median of $f$ (assume it exists), defined by simultaneously
satisfying 
\[
\mu\{ x\in \Omega_n : f(x)\leq M_f\} \geq 1/2
\]
and 
\[
\mu\{ x\in \Omega_n  : f(x)\geq M_f\} \geq 1/2.
\]
The ramifications of definition (\ref{concentrationfunction}) and property (\ref{normalconcentration}) can be seen in considering the
set $A:= \{x\in \Omega_n : f(x)=M_f\}$. 
One may show that \cite{asymp}
\[
\mu(A_{\epsilon})
\geq 1-2\alpha(\Omega_n,\epsilon).
\]
In particular, if $(\Omega_n,\mu)$ has normal concentration (\ref{normalconcentration}), we have
\begin{equation}\label{nconcentrate}
\mu(A_{\epsilon}) \geq 1-2c_1\exp(-c_2\epsilon^2n).
\end{equation}
Now if we assume that $\omega_f(\epsilon)\leq \delta$ is small, then
(\ref{nconcentrate})
states that the values of $f$ are very close
to $M_f$ on almost all of the space. Indeed, given $x\in A_{\epsilon}$,
there exists $y\in A$ such that $d(x,y)\leq \epsilon$. Then 
$|f(x)-f(y)| \leq \omega_f(\epsilon)$, and since $y\in A$, 
$|f(x)-M_f| \leq \omega_f(\epsilon)$. 
We conclude the distribution of $f$ satisfies the {\it concentration of measure phenomenon}:
\begin{eqnarray*}
\mu\{ x\in \Omega_n:|f(x)-M_f| \leq \omega_f(\epsilon)\}&\geq&
\mu(A_{\epsilon}) \\
&\geq&
1- 2\alpha(\Omega_n,\epsilon)\\
&\geq& 1-2c_1\exp(-c_2\epsilon^2n).
\end{eqnarray*}
(We note that one may interchange the median for the mean if $f$ is 
Lipschitz.)

We have passed  from an initial geometric formulation of the concentration of measure phenomenon involving $\epsilon$-neighborhoods of subsets of measure at least 1/2, to a purely probabilistic one involving tail bounds of Lipschitz functions about their mean.
The geometric formulation has its roots in the study of isoperimetric inequalities or $\R^n$: among bodies of a given volume, Euclidean balls have the least surface area \cite{ball}. Since the idea of comparing volumes of a set and its neighborhoods makes sense in any metric probability space $(\Omega,d)$, one may pose the following more abstract question: for which sets $A$ of measure $\alpha$ do the $\epsilon$-neighborhoods $A_{\epsilon}$ have the smallest measure? This is known as the isoperimetric problem. An initial answer to this question, which relates to present day understanding of the concentration of measure phenomenon, was given by P. L\'evy. He considered the isoperimetric problem on the sphere $S^{n-1}$, equipped with the Haar (rotation-invariant) probability measure $\sigma$, and the geodesic distance. He proved that in this case it was spherical caps that solve the isoperimetric problem. Thus, if $A$ is a subset of the sphere of the same measure as a spherical cap, then its $\epsilon$-neighborhood is at least as large in measure as that of the $\epsilon$-neighborhood of the cap.  If the measure of the set $A$ is 1/2, then upon comparing $A_{\epsilon}$ with the $\epsilon$-neighborhood of a hemisphere (spherical cap of measure 1/2), one obtains the bound $\sigma(A_{\epsilon}) \geq 1-e^{-n\epsilon^2/2}$ (a purely computational fact, though rather difficult \cite{milman, ball}). Thus almost the entire sphere (measure-theoretically) lies within a distance $\epsilon$ of $A$. 

As Talagrand highlighted in \cite{talagrand},  the general probabilistic notion of the concentration of measure phenomenon describes the tail behavior of random variables which depend in a ``smooth'' way on many (perhaps independent) random variables (in a way that generalizes the central limit theorem).
This was studied in great detail by Talagrand \cite{talagrand} on an abstract product probability space $(\Omega^n, \p)$, with $\p$ a product measure. There is a natural metric on any product space, namely the Hamming distance
\[
d_H(x,y) = \text{card}\{i\leq n: x_i\neq y_i\},
\]
so it is reasonable to begin with this to study concentration phenomena on product probability spaces. 
In this setting Talagrand proved
\begin{equation}\label{Hamming}
\int_{\Omega^n} e^{td_H(x,A)} \, d\p(x) \leq \frac{1}{\p(A)} e^{t^2n/4},
\end{equation}
which, by an application of Markov's inequality, gives bounds on the probability that a point lies outside or inside an $\epsilon$-neighborhood of $A$.
Remarkably, Talagrand was able to prove a similar inequality for the Euclidean distance
 $d$ on $\R^{2n}$ (actually
a more general ``combinatorial'' distance), proving
\begin{equation}\label{euclidean}
\int_{\R^n} e^{\frac{1}{4} d^2(A,x)} \, d\p(x) \leq \frac{1}{\p(A)}.
\end{equation}
For both (\ref{Hamming}) and (\ref{euclidean}) Talagrand's argument exploited the inductive nature (in $n$) of the product measure along with the inductive nature of the particular distance. 

However, such an approach seems particular to the product distribution. Other methods have been developed to deal with non-product distributions, many of which depend on bounding the dependencies of the variables in terms of various types of mixing coefficients \cite{Marton,marton3,marton2, Samson}. Other research on the concentration of non-product distributions can be found for example in \cite{Kontor,peres}. 

In this paper   
we develop a method of induction to deal with non-product distributions to study a simplified version of the concentration of measure phenomenon that retains a geometric appeal. We consider the concentration of the  metric-induced function, $x\mapsto d(x,y)$, $y\in \Omega$, i.e., the distribution of distances about a fixed $y$. That is, we consider the concentration of the probability space itself.
Our results are specialized to the Boolean cube $I_n:= \{0,1\}^n$, equipped with Hamming distance $d_H$, as well as Euclidean distance. The Boolean cube is a standard simple model on which concentration phenomena have been studied \cite{talagrand,spencer}; often phenomena in this setting have implied similar ones on more general spaces. This was the case, for example, in \cite{talagrand}. The space is finite, so measure-theoretic subtleties are of no concern, but it has enough structure to provide non-trivial problems.
A similar rich source of questions regarding concentration inequalities comes from the theory of random graphs \cite{spencer}.

For a given $y\in I_n$, the Hamming distance function, $x\mapsto d_H(x,y)$, is Lipschitz in the variable $x$.  
In fact it is a sum of bounded random variables,
\[
d_H(x,y) = \sum_{i=1}^n |x_i-y_i|.
\]
Therefore the idea that $x\mapsto d_H(x,y)$ concentrates about its mean is not a novel notion, at least in the situation when $x$ has independent components. 

But suppose $x$ is not independent in each component.
An initial observation is that the concentration of $d_H$ fails if the measure $\mu$
``separates mass''. For example take $\mu = \frac{1}{2}\delta_{(0,\dots,0)}
+\frac{1}{2}\delta_{(1,\dots,1)}$, which separates mass evenly between the two points $(0,\dots, 0)$
and $(1,\dots, 1)$. For any $y\in I_n$ we have 
$\E_{\mu} d_H(x,y) = n/2$; however the single point
$(0,\dots, 0)$ has measure 1/2, and therefore 
$d_H(x,0)$ is zero on a subset of measure 1/2,
at all dimensions $n\geq 1$. Thus $d_H(x,0)$ 
cannot concentrate about $n/2$. 

We note that in this case, $\mu(x)=0$ for a majority
of the points on $I_n$. However, supposing that $\mu(x)>0$ for all $x\in I_n$, one
can ask whether concentration can hold under this further constraint. If we fix $\epsilon>0$ 
sufficiently small and consider the distribution 
\[
\mu= \left(\frac{1}{2}-\epsilon\right)\delta_{(0,\dots,0)}+\frac{1}{2}
\delta_{(1,\dots,1)}+\frac{\epsilon}{2^n-2}\sum_{x\neq \overline{0},\overline{1}}
\delta_{x}
\]
where $\overline{0}=(0,\dots, 0)$, $\overline{1}=(1,\dots,1)$.
This distribution still separates mass between $\overline{0}$
and $\overline{1}$, but now an $\epsilon$ amount is
evenly distributed to the remaining points. Thus we have
$\mu(x)>0$ for all $x\in I_n$. Given any $y\in I_n$ with
$k$ ones among its components, we see that the mean value
of $d_H(x,y)$ with respect to $\mu$ is 
\begin{eqnarray*}
\E_{\mu}d_H(x,y)&=& \left(\frac{1}{2}-\epsilon\right)k +\frac{1}{2}(n-k)
+\frac{\epsilon}{2^n-2}\sum_{x\neq \overline{0},\overline{1}}
d_H(x,y) \\
&=&\frac{n}{2}-k\epsilon + \frac{\epsilon}{2^n-2}\sum_{x\neq \overline{0},\overline{1}}
d_H(x,y).
\end{eqnarray*}
In particular, if $y=\overline{0}$ then its mean is rather close to 
$n/2$, but $d_H(x,\overline{0})$ is zero on a subset of measure 
$1/2-\epsilon$. Thus $d_H$ cannot concentrate about its mean
with respect to $\mu$.

Now instead of demanding that $d_H(x,y)$ concentrate about its mean
for all $y\in I_n$, we rather consider a fixed $y\in I_n$. The above
examples show that  the quality of concentration can depend
on the choice of $y$. In section \ref{PAsection} we will in fact see that
this is very much the case for {\it all} dependent distributions $\mu_n$ on
the Boolean cube. It is only with respect to (close to) independent distributions
that one obtains a uniform upper bound of the form (usually going by the name of ``Hoeffding's inequality'' \cite{CM}),
\begin{equation}\label{Hoeff}
\int_{I_n} e^{t(d_H(x,y) - \E_{\mu_n} d_H(x,y))} \, d\mu_n(x) \leq e^{nt^2/2}
\end{equation}
for all $y\in I_n$. The key difference in the dependent case is that
the conditional distributions 
\[
\mu(x_k \mid x_1,\dots, x_{k-1}), \;\;\; 2\leq k \leq n,
\]
are random quantities, in that they depend on the random variables $x_1,\dots, x_{k-1}$. However, we can consider their mean values
\[
\E_{\mu_{k-1}} \mu(x_k \mid x_1,\dots, x_{k-1}) = \mu^{(k)}(x_k),
\]
where $\mu^{(k)}$, $1\leq k\leq n$, denote the one-dimensional marginal
distributions of $\mu_n$ ( distribution of the variable $x_k$, $1\leq k\leq n$), and $\mu_k$ denotes the joint distribution of $x_1,\dots, x_k$. 
The above conditional distributions being random quantities, it is in fact useful to quantify {\it how far the full measure
$\mu_n$ is from being independent} by defining the (mean
zero) variables $\epsilon_{x',y_k}$, $x'=(x_1,\dots, x_k)$, given
by 
\[
\mu(x_k=y_k \mid x_1,\dots, x_{k-1}) = \mu^{(k)}(x_k=y_k) +\epsilon_{x',y_k},
\]
with the left side containing a conditional probability that $x_k=y_k$, and the right side the marginal probability of the same event. These essentially measure the deviation from independence on the introduction of the $k^{th}$ component $x_k$ to the first $k-1$ components.
When we change $y_k$, the sign of $\epsilon_{x',y_k}$ changes.
It turns out that this subtle yet simple property will dictate the sign of an
``error-term'' in a bound we derive in {\bf Theorem \ref{errorbound}}: 
\begin{equation}\label{boundintro}
\int_{I_n} e^{t(d_H(x,y) - \E_{\mu_n} d_H(x,y))} \, d\mu_n(x) \leq e^{nt^2/2}+(\text{error}), \;\;\;\; t>0.
\end{equation}
The factor of $e^{nt^2/2}$ is the same as that in the independent case, studied in (\ref{Hoeff}). Further, the ``error'' term we derive is inductive, and of the form
\[
\text{error} = \sum_{k=1}^{n-1} e^{(n-k-1)t^2/2}E_{k,x',y',t}
\]
where 
\[
E_{k,x',y',t} = e^{-t\mu^{(k+1)}(x_{k+1}\neq y_{k+1})}(1-e^t) \int_{I_{k}}  e^{t(d_H(x',y')-\E d_H(x,y))}\epsilon_{x',y_{k+1}}\, d\mu_{k}(x').
  \]
We see then that the {\it correlation} between $e^{t(d_H(x',y')-\E d_H(x,y))}$ and $\epsilon_{x',y_{k+1}}$ dictates the sign of this error-term. In particular, if we have {\it positive correlation} with respect to $y_{k+1}$, then we must have {\it negative correlation} with respect to $1-y_{k+1}$.
And therefore the quality of concentration will depend on the choice of $y$. A simple consequence of this is that, for any distribution on the Boolean cube $I_n$, there will always exists a $y\in I_n$ for which the function $x\mapsto d_H(x,y)$ (perhaps trivially) concentrates about its mean.
In the worst case, when for each $k\geq 1$ we have $E_{k,x',y',t}>0$, there is a rather elegant product formula for (\ref{boundintro}). It takes a form which in this paper is {\bf Theorem \ref{smallvariancebound}},
\[
\int_{I_n} e^{t(d_H(x,y)-\E d_H(x,y))} \, d\mu(x) \leq e^{t^2/2} \prod_{j=2}^{n} (b_j +e^{t^2/2}), \;\;\;\; t>0,
\]
where $b_j := e^{-t\mu^{(j)}(x_j\neq y_j)}(e^t-1)\lVert \epsilon_{x',y_j}\rVert_{\infty}.$ In particular it highlights in a simple fashion how the variance of the conditional distributions, or equivalently the size of $\epsilon_{x',y_k}$, affects the concentration bound.

An outline of this paper is as follows. In section 2, broadly, we develop our induction method and prove an error bound for concentration (that is a bound on the exponential moment generating function about $e^{nt^2/2}$). In section \ref{positivecorrelation} we introduce a {\it positive correlation} condition on the coordinates of the point $y\in I_n$ (fixed in $x\mapsto d_H(x,y)$) which will dictate the sign of the error term, and thus the quality of concentration. We then show, in section \ref{smallvariance}, how the error bound simplifies elegantly in the worst case, deriving what we call the small-variance bound.
Finally, in section \ref{Lipschitzuniform} we demonstrate how the small- variance bound generalizes to the entire class of Lipschitz functions on the Boolean cube.

\section{The Concentration of Hamming Distance}

\subsection{Preliminaries on the Boolean Cube}

Let $x_1,\dots,x_n$ denote 0-1 valued random variables with joint distribution
$\mu_n$. Thus the vector $(x_1,\dots, x_n)$ may be viewed as an element
in the Boolean cube $I_n = \{0,1\}^n$, and the distribution $\mu_n$ a probability
measure on $I_n$.  The distribution of each variable $x_i$, $i=1,\dots, n$, is denoted
by $\mu^{(i)}$, these are the one-dimensional marginal distributions of
$\mu_n$.
For any fixed $x_1,\dots, x_{k-1}$, $2\leq k\leq n$, we have
\[
\mu_k(x_1,\dots, x_{k-1},0) + \mu_k(x_1,\dots, x_{k-1},1)
= \mu_{k-1} (x_1,\dots, x_{k-1}).
\]
Thus there exist $0\leq c_{x_1,\dots, x_{k-1}} \leq 1$ such that 
\begin{equation}\label{marginal}
\mu_k(x_1,\dots,x_{k-1}, 0)= c_{x_1,\dots,x_{k-1}}\mu_{k-1}(x_1,\dots, x_{k-1}),
\end{equation}
and 
\begin{equation}\label{marginal2}
\mu_k(x_1,\dots, x_{k-1},1)= (1-c_{x_1,\dots,x_{k-1}})\mu_{k-1}(x_1,\dots, x_{k-1}).
\end{equation}
That is, 
\[
c_{x_1,\dots, x_{k-1}} = \mu(x_k=0 | x_1,\dots, x_{k-1}).
\]
Moreover, we see from (\ref{marginal}) (or from the fact that $c_{x_1,\dots,x_{k-1}}=\mu(x_k=0 | x_1,\dots,x_{k-1})$) that
\[
\E_{\mu_{k-1}} c_{x_1,\dots,x_{k-1}}= \mu^{(k)}(0) \;\;\;\; \forall k=2,\dots n.
\]
Note that all the expectations above and below are with respect to the variable $x\in I_n$, while $y\in I_n$ is fixed. 
Now fix $y\in I_{n}$. Noting that $d_H(x,y) = \sum_{i=1}^n |x_i-y_i|$, we have 
\begin{equation}\label{averagehamming}
\E_{\mu_n} d_H(x,y) = \sum_{i=1}^n \E_{\mu^{(i)}} |x_i-y_i| = \sum_{i=1}^n \mu^{(i)}(|1-y_i|).
\end{equation}
Note that without abuse of notation, we will write $\mu^{(i)}(x_i=|1-y_i|) \equiv \mu^{(i)}(|1-y_i|)$, i.e. the variable $x$ will be omitted from the arguments of $\mu$ when possible.
We will also use the notation, $\mu^{(i)}(|1-y_i|) = \mu^{(i)}( x_i\neq y_i).$
We will moreover need a subtle yet powerful lemma known as Hoeffding's lemma  \cite{CM, TT}:
\begin{lemma}\label{hoeffding}
Let $f$ be integrable with respect to a probability distribution $\p$, and assume $\E[f]=0$. 
Then 
\[
\E[e^{tf}] \leq \frac{e^{t \lVert f \rVert_{\infty}} + e^{-t \lVert f \rVert_{\infty}}}{2} \leq 
e^{t^2\lVert f\rVert_{\infty}^2/2}.
\]
\end{lemma}

\subsection{A Classification of the Concentration of the Hamming Distance}\label{PAsection}
It is known that for independent measures on the Boolean cube, standard concentration inequalities for Lipschitz functions on the measure space hold, but these do not extend in known ways to more general Boolean cube measures. Here we show that in fact such results can be extended, using the present notion of concentration.

For a given $y\in I_n$ the Hamming distance $d_H(x,y)$ is Lipschitz (with Lipschitz constant one). At a more basic level, it is actually a sum of bounded random variables:
\[
d_H(x,y) = \sum_{i=1}^n |x_i-y_i|.
\]
Thus when the vector $(x_1,\dots,x_n)$ has independent components, concentration is well known, and usually goes by the name of Hoeffding's inequality (as above) \cite{CM}:  
\begin{equation}\label{independentcase}
\int_{I_n}e^{t(d_H(x,y) -\E d_H(x,y))} \, d\mu(x)\leq e^{nt^2/2}
\end{equation}
for any $\mu$ that decomposes as a product of its one-dimensional marginals.
Thus by Markov's inequality (setting $t=c/n$ and using the symmetry in $t$),
\begin{equation}\label{independenttailbound}
\mu\{ x\in I_n :|d_H(x,y)-\E d_H(x,y)| \geq c \} \leq 2e^{-c^2/2n}.
\end{equation}
In particular, the ``observable diameter" of the Boolean cube is of order $\sqrt{n}$. That is, for any fixed $y\in I_n$ and $n$ large, $d_H(x,y)$ lies within a distance of $\sqrt{n}$ of its mean on a subset of measure converging to one in a sub-Gaussian fashion (set $c=\sqrt{n}$). The bound in (\ref{independentcase}) is quite optimal in that it produces sub-Gaussian tail bounds (\ref{independenttailbound}). But much is unknown as to whether, or when, such bounds hold for dependent distributions, i.e. dependent vectors $(x_1,\dots, x_n)$ and their distributions. 

In this section, for a given $y\in I_n$, we classify the concentration of the Hamming distance,
$d_H(x,y)$, about its mean, with respect to {\it arbitrary} distributions $\mu_n$ on the Boolean cube. We will see
that the quality of concentration depends on essentially two factors. These are
the point $y\in I_n$ chosen and the variance of the conditional distributions
$\mu(x_k \mid x_1,\dots, x_{k-1})$, $2\leq k\leq n$. (These variances can be thought of as how far $\mu$ is from
being independent.)

\subsubsection{An Inductive Error Bound Detailing Concentration}

Fix $2\leq k\leq n$, $y_k\in \{0,1\}$,
and consider $\mu(x_k=y_k \mid x_1,\dots, x_{k-1})$. The mean value
of this random variable is 
\begin{equation}
\E_{\mu_{k-1}} \mu(x_k=y_k \mid x_1,\dots, x_{k-1})
=\mu^{(k)}(x_k=y_k).
\end{equation}
 Thus 
\begin{equation}\label{conditionalmeasure}
\mu(x_k=y_k \mid x_1,\dots, x_{k-1}) = \mu^{(k)}(x_k=y_k) +\epsilon_{x_1,\dots, x_{k-1},y_k}
\end{equation}
and 
\begin{equation}
\mu(x_k\neq y_k \mid x_1,\dots, x_{k-1}) = \mu^{(k)}(x_k\neq y_k) +\epsilon_{x_1,\dots, x_{k-1},1-y_k}
\end{equation}
where $\E_{\mu_{k-1}} \epsilon_{x_1,\dots, x_{k-1},y_k}=0$. 
Since $1-\mu(x_k=y_k \mid x_1,\dots, x_{k-1}) =\mu(x_k\neq y_k\mid x_1,\dots, x_{k-1})$,
and $1-\mu(x_k=y_k )=\mu(x_k\neq y_k)$, we have
the important relation, 
\begin{equation}
\epsilon_{x_1,\dots, x_{k-1},y_{k}} = -\epsilon_{x_1,\dots, x_{k-1},{1-y_{k}}}.
\end{equation}
The variables $\epsilon_{x_1,\dots, x_{k-1},y_k}$ define the
spread of the distribution of $\mu(x_k = y_k \mid x_1,\dots, x_{k-1})$
over all $(x_1,\dots, x_{k-1})\in I_{k-1}$, and
therefore the variance of this distribution, and so how far the distribution
is from being independent of $(x_1,\dots, x_{k-1})$.

For a given $2\leq k \leq n$,
define $x' = (x_1,\dots, x_k)$, $y'=(y_1,\dots, y_k)$, and (again for fixed $y$)
\begin{equation}
a_{k-1}(x',y',t) := e^{t(d(x',y') -\E_{\mu_{k-1}} d(x',y'))}.
\end{equation}

\begin{thm}\label{errorbound}
For any probability distribution $\mu_n$ on the Boolean cube, with $y$ fixed and expectations over $x$,
\begin{equation}\label{errorterm}
\int_{I_n} e^{t(d_H(x,y) -\E d_H(x,y))} \, d\mu_n(x)\leq
e^{nt^2/2} +\sum_{k=1}^{n-1}e^{(n-k-1)t^2/2}E_{k,x',y',t}
\end{equation}
where
\[
 E_{k,x',y',t} := e^{-t\mu^{(k+1)}(x_{k+1}\neq y_{k+1})}(1-e^t) \int_{I_{k}}  a_{k}(x',y',t)\epsilon_{x',y_{k+1}}\, d\mu_{k}(x').
  \] 
  \end{thm}

\begin{Proof}
We proceed by induction. The case of $n=1$ is Hoeffding's lemma, so assume
the result in dimension $n-1$.
Let $y=(y_1,\dots,y_n)\in I_n$ be fixed, and consider the following decomposition,
\begin{equation}\label{decomp}
\int_{I_n} e^{t(d_H(x,y) -\E_{\mu_n} d_H(x,y))} \, d\mu_n(x) =\left(
\int_{I_n^{y_n}} + \int_{I_n^{1-y_n}}\right) e^{t(d_H(x,y) -\E_{\mu_n} d_H(x,y))} \, d\mu_n(x) ,
\end{equation}
where $I_n^{y_n}$ denotes all $x=(x_1,\dots, x_n)\in I_n$ such that $x_n=y_n$,
and $I_n^{1-y_n}$ denotes all $x\in I_n$ such that $x_n\neq y_n$.
Let $x' = (x_1,\dots,x_{n-1})$ and $a_{n-1}(t,x',y') = e^{t(d_H(x',y') - \E_{\mu_{n-1}}d_H(x',y'))}$ as above.
We have
\begin{eqnarray*}
&&\int_{I_n^{y_n}} e^{t(d_H(x,y) -\E_{\mu_n} d_H(x,y))} \, d\mu_n(x) =
\sum_{x\in I_n^{y_n}}  e^{t(d_H(x,y) -\E_{\mu_n} d_H(x,y))} \mu_n(x) \\
&=& \sum_{x'\in I_{n-1}}  a_{n-1}(x',y',t) e^{t(|y_n-y_n| -\E_{\mu^{(n)}}|x_n-y_n|)}
\mu_n(x',y_n) \\ 
&=&  \sum_{x'\in I_{n-1}}  a_{n-1}(x',y',t) e^{t(-\mu^{(n)}(x_n\neq y_n))}\mu_n(x',y_n)\\
&=& \sum_{x'\in I_{n-1}}  a_{n-1}(x',y',t) e^{-t\mu^{(n)}(x_n\neq y_n)}
\mu_n(x_n=y_n\mid x')\mu_{n-1}(x') 
\end{eqnarray*}
Using that  $\mu_n(x_n=y_n\mid x') = \mu^{(n)}(x_n=y_n) + \epsilon_{x',y_n}$,
for all $x'\in I_{n-1}$,
we may re-write this last expression as
\[
 \sum_{x'\in I_{n-1}}  a_{n-1}(x',y',t) e^{-t\mu^{(n)}(x_n\neq y_n)}\mu^{(n)}(x_n=y_n) \mu_{n-1}(x')+ 
 \]
 \[
 \sum_{x'\in I_{n-1}}  a_{n-1}(x',y',t)e^{-t\mu^{(n)}(x_n\neq y_n)} \epsilon_{x',y_n}\mu_{n-1}(x').
 \]
Considering now the integral over $I_n^{1-y_n}$, we have,
\begin{eqnarray*}
&&\int_{I_n^{1-y_n}} e^{t(d_H(x,y) -\E_{\mu_n} d_H(x,y))} \, d\mu_n(x) =
\sum_{x\in I_n^{1-y_n}}  e^{t(d_H(x,y) -\E_{\mu} d_H(x,y))} \mu_n(x) \\
&=& \sum_{x'\in I_{n-1}}  a_{n-1}(x',y',t) e^{t(|(1-y_n)-y_n| -\E_{\mu^{(n)}}|x_n-y_n|)}
\mu_n(x',1-y_n) \\ 
&=&  \sum_{x'\in I_{n-1}}  a_{n-1}(x',y',t) e^{t(1-\mu^{(n)}(x_n\neq y_n))} \mu_n(x',1-y_n)\\
&=& \sum_{x'\in I_{n-1}}  a_{n-1}(x',y',t) e^{t(1-\mu^{(n)}(x_n\neq y_n))}
\mu_n(x_n\neq y_n\mid x')\mu_{n-1}(x')
\end{eqnarray*}
Note that the third equality above uses the fact that $|(1-y_n)-y_n| = |1-2y_n|=1$ for $y_n=0,1$.
Since $\mu_n(x_n=y_n\mid x') = \mu^{(n)}(x_n=y_n) + \epsilon_{x',y_n}$,
we see that $\mu_n(x_n\neq y_n\mid x') = 1-\mu^{(n)}(x_n=y_n) - \epsilon_{x',y_n}
=\mu^{(n)}(x_n \neq y_n) - \epsilon_{x',y_n}$. Thus the last expression
in the above string of equalities may be re-written as 
\[
  \sum_{x'\in I_{n-1}}  a_{n-1}(x',y',t) e^{t(1-\mu^{(n)}(x_n\neq y_n))} \mu^{(n)}(x_n \neq y_n) \mu_{n-1}(x') -
  \]
  \[
   \sum_{x'\in I_{n-1}}  a_{n-1}(x',y',t) e^{t(1-\mu^{(n)}(x_n\neq y_n))} \epsilon_{x',y_n} \mu_{n-1}(x') .
\]
Thus, using these results in (\ref{decomp}), we find that $\int_{I_n}e^{t(d_H(x,y)-\E_{\mu_n}d_H(x,y))}\, d\mu_n(x)$ can be written as 
\[
 \sum_{x'\in I_{n-1}}  a_{n-1}(x',y',t) e^{-t\mu^{(n)}(x_n\neq y_n)}(\mu^{(n)}(x_n=y_n) +e^t\mu^{(n)}(x_n\neq y_n))\mu_{n-1}(x') +
 \]
 \[
  \sum_{x'\in I_{n-1}}  a_{n-1}(x',y',t) e^{-t\mu^{(n)}(x_n\neq y_n)}\epsilon_{x',y_n}(1-e^t)\mu_{n-1}(x').
  \]
Note that  $e^{-t\mu^{(n)}(x_n\neq y_n)}(\mu^{(n)}(x_n=y_n) +e^t\mu^{(n)}(x_n\neq y_n))$
can be re-written as an integral, and then bounded above using the one-dimensional version of Hoeffding's lemma \ref{hoeffding}, i.e. using the identity
\[
\int_{I_1} e^{t(d_H(x_n,y_n) -\E_{\mu^{(n)}}d_H(x_n,y_n))}\, d\mu^{(n)}(x_n)
\leq e^{t^2/2}
\]
(note that this is a one-dimensional integral in $x_n$). This uses the fact that $ 
e^{-t E_{\mu^{(n)}}d_H(x_n,y_n)} = e^{-t\mu^{(n)}(x_n\neq y_n)}$, and also that  $\mu^{(n)}(x_n=y_n)+ e^t\mu^{(n)}(x_n\neq y_n) 
= \int_{I_1} e^{t d_H(x_n,y_n)} d\mu^{(n)}(x_n)$.
Thus $\int_{I_n} e^{t(d_H(x,y) -\E_{\mu} d_H(x,y))} \, d\mu(x) $ is bounded above
by (changing back to integral notation),
\[
e^{t^2/2} \int_{I_{n-1}} a_{n-1}(x',y',t)\, d\mu_{n-1}(x')+E_{n-1,x',y',t}
\]
where 
\[
E_{n-1,x',y',t}=e^{-t\mu^{(n)}(x_n\neq y_n)}(1-e^t) \int_{I_{n-1}}  a_{n-1}(x',y',t) \epsilon_{x',y_n} \, d\mu_{n-1}(x').
  \]

Using the induction step (\ref{errorterm}) in dimension $n-1$, we have
\begin{eqnarray*}
&&e^{t^2/2}\int_{I_{n-1}} a_{n-1}(x',y',t) \, d\mu_{n-1}(x') + E_{n-1,x',y',t}\\
&\leq&e^{nt^2/2}+\left( \sum_{k=1}^{n-2}e^{(n-k-1)t^2/2}E_{k,x',y',t}\right)
+ E_{n-1,x',y',t} \\
&=&e^{nt^2/2}+\sum_{k=1}^{n-1}e^{(n-k-1)t^2/2}E_{k,x',y',t}.
\end{eqnarray*}
The result follows. 
\end{Proof}

\subsubsection{A Positive Correlation Condition for Concentration}\label{positivecorrelation}

The error bound given in Theorem \ref{errorbound} (last term in (\ref{errorterm}) seems 
hard to unravel, but there is a key observation that allows it to be better understood:

\begin{prop}\label{positiveprop}
 For a fixed $1\leq k\leq n-1$, the sign
of the error-term, $E_{k,x',y',t}$, depends on the choice of $y_{k+1}$.
\end{prop}

\begin{Proof}
Fix $1\leq k\leq n-1$. Then 
\[
E_{k,x',y',t} = e^{-t\mu^{(k+1)}(x_{k+1}\neq y_{k+1})}(1-e^t) \int_{I_{k}}  a_{k}(x',y',t)\epsilon_{x',y_{k+1}}\, d\mu_{k}(x') \leq 0
  \] 
if and only if,
\begin{equation}\label{positiveterm1}
 \int_{I_{k}}  a_{k}(x',y',t) \epsilon_{x',y_{k+1}}\, d\mu_{n-1}(x') \geq 0
 \end{equation}
 since $1-e^t<0$ for all $t>0$. Since $\epsilon_{x',y_{k+1}} = -
 \epsilon_{x',1-y_{k+1}}$, for all $x'\in I_{k}$, we see that (\ref{positiveterm1}) is equivalent to
 \[
  \int_{I_{k}}  a_{k}(x',y',t) \epsilon_{x',1-y_{k+1}} \, d\mu_{n-1}(x') \leq 0.
  \]
Thus for fixed $k$ the sign of the error term $E_{k,x',y',t}$ is opposite for the two 
choices of $y_{k+1}$.
 \end{Proof}

In light of Proposition \ref{positiveprop}, for a given $2\leq k\leq n$,
we make the following definition:

\begin{defn}
For a given $2\leq k\leq n$,
we say that $y_k\in\{0,1\}$ satisfies the {\it positive correlation condition} if
\[
\int_{I_{k-1}} a_{k-1}(x',y',t) \mu(x_k=y_k\mid x')\,  d\mu_{k-1}(x') \geq 
\]
\begin{equation}\label{PA}
\int_{I_{k-1}} a_{k-1}(x',y',t)\, d\mu_{k-1}(x') \int_{I_{k-1}} \mu(x_k=y_k \mid x') \, d\mu(x') 
\end{equation}
Since $ \int_{I_{k-1}} \mu(x_k=y_k \mid x') \, d\mu(x') =\mu^{(k)}(x_k=y_k)$,
we may equivalently write (\ref{PA}) as,
\begin{equation}
\int_{I_{k-1}} a_{k-1}(x',y',t) \epsilon_{x',y_k}\,  d\mu_{k-1}(x') \geq  0.
\end{equation}
\end{defn}

We say that this is a positive correlation condition, because (\ref{PA}) states
that the variables $X := a_{k-1}(x',y',t)$ and $Y:= \mu(x_k=y_k \mid x')$
are positively correlated.

Since $\mu(x_k=y_k \mid x') = 1-\mu(x_k\neq y_k \mid x')$, (\ref{PA})
implies the following {\it negative correlation condition} for $1-y_k$:

\begin{lemma}
If $y_k$ satisfies the positive correlation condition (\ref{PA}), then $1-y_k$ satisfies the following negative correlation condition:
\[
\int_{I_{k-1}} a_{k-1}(x',y',t) \mu(x_k\neq y_k\mid x')\,  d\mu_{k-1}(x') \leq
\]
\begin{equation}\label{NAcondition}
 \int_{I_{k-1}} a_{k-1}(x',y',t)\, d\mu_{k-1}(x') \int_{I_{k-1}} \mu(x_k\neq y_k \mid x') \, d\mu(x').
\end{equation}
Since  $\int_{I_{k-1}} \mu(x_k\neq y_k \mid x') \, d\mu(x')=\mu^{(k)}(x_k\neq y_k)$,
the negative correlation condition (\ref{NAcondition}) may be equivalently written, 
\begin{equation}
\int_{I_{k-1}} a_{k-1}(x',y',t) \epsilon_{x',1-y_k}\,  d\mu_{k-1}(x') \leq  0.
\end{equation}
\end{lemma}

For example, if $y_n=1$ and we have the positive correlation
condition
\[
\int_{I_{n-1}} a_{n-1}(x',y',t) \epsilon_{x',1}\,  d\mu_{n-1}(x') \geq  0
\]
then necessarily we have the negative correlation condition for 
$y_n=0$,
\[
\int_{I_{n-1}} a_{n-1}(x',y',t) \epsilon_{x',0}\,  d\mu_{n-1}(x') \leq  0
\]

Proposition \ref{positiveprop} shows that the positive and negative correlation conditions above dictate the quality of concentration for a given $y\in I_n$. In particular,
the best concentration occurs when $y$ is such that the positive correlation
condition (\ref{PA}) holds for all $k$, with the bound from above the same as in the case that the variables are independent: 

\begin{thm}\label{PAtheorem}
Suppose $y=(y_k)_{k=1}^n\in I_n$ is such that the positive correlation condition (\ref{PA}) holds for all $k\geq 2$. Then
\[
\int_{I_n} e^{t(d_H(x,y) -\E_{\mu} d_H(x,y))} \, d\mu(x) \leq e^{nt^2/2}.
\]
\end{thm}

\begin{Proof}
It suffices to show that the error terms $E_{k,x',y',t}$ are at most
zero for all $k\geq 2$. We have 
\[
E_{k,x',y',t} = e^{-t\mu^{(k+1)}(x_{k+1}\neq y_{k+1})}(1-e^t) \int_{I_{k}}  a_{k}(x',y',t) \epsilon_{x',y_{k+1}}\, d\mu_{k}(x') \leq 0
\]
if and only if,
\begin{equation}\label{positiveterm}
 \int_{I_{k}}  a_{k}(x',y',t) \epsilon_{x',y_{k+1}}\, d\mu_{k}(x') \geq 0
 \end{equation}
 since $1-e^t<0$ for all $t>0$. Since (\ref{positiveterm}) is exactly the
 positive correlation condition, the result follows.
\end{Proof}

If the error terms $E_{k,x',y',t}$ in Theorem \ref{errorbound} are proportional at all levels, in the sense that
\[
e^{(n-k-1)t^2}|E_{k,y',x',t}|,\;\;\; k=1,\dots, n-1,
\]
is constant in $k$, then the sum 
\begin{equation}\label{sumerror}
\sum_{k=1}^{n-1} e^{(n-k-1)t^2/2}E_{k,x',y',t}
\end{equation}
will be approximately bounded above by zero for a large number of $y\in I_n$. Specifically it will hold for those $y$ for which more coordinates $y_k$ satisfy the positive correlation condition (\ref{PA}) than the negative correlation condition. We state this more precisely as a theorem.

\begin{thm}
 Suppose $\mu$ is a probability distribution on the Boolean cube satisfying 
 \[
 \mu^{(i)}(0) = 1/2, \;\;\; i=1,\dots, n.
 \]
 Moreover suppose that 
 \[
e^{(n-k-1)t^2}|E_{k,y',x',t}|
\]
is constant in $k$, $1\leq k\leq n-1$. Then there exist $2^{n-\lceil\frac{n-1}{2}\rceil}\binom{n}{\lceil\frac{n-1}{2}\rceil}$ vectors $y\in I_n$, with respect to which the function $x\mapsto d_H(x,y)$ satisfies
\[
\int_{I_n} e^{t(d_H(x,y) -\E d_H(x,y))} \, d\mu(x) \leq e^{nt^2/2}.
\]
\end{thm}

\begin{Proof}
Denote $E_{k,y_k} = E_{k,x',y',t}$, so
\[
E_{k,y_k} = e^{-t\mu^{(k+1)}(x_{k+1}\neq y_{k+1})}(1-e^t) \int_{I_{k}}  a_{k}(x',y',t)\epsilon_{x',y_{k+1}}\, d\mu_{k}(x').
  \] 
We already know that $\epsilon_{x',y_{k+1}}=-
\epsilon_{x',1-y_{k+1}}$,
so $E_{k,y_k}$ is approximately
$-E_{k,1-y_k}$. 
But with the added assumption
that $\mu^{(i)}(0)=1/2$
for each $i=1,\dots, n$, we have exact equality:
\[
E_{k,y_k} = - E_{k,1-y_k}
\]
for each $k=1,\dots, n-1$. And since
 \[
e^{(n-k-1)t^2}|E_{k,y',x',t}|
\]
is constant in $k$, $1\leq k\leq n-1$, the function $x\mapsto d_H(x,y)$ will satisfy 
\[
\int_{I_n} e^{t(d_H(x,y) -\E d_H(x,y))} \, d\mu(x) \leq e^{nt^2/2},
\]
if $E_{k,y_{k+1}}\leq 0$ for at least half the coordinates of $y$ (since the negative error terms will cancel the positive error terms in (\ref{sumerror})). So if we choose $\binom{n}{\lceil\frac{n-1}{2}\rceil}$ coordinates $y_k$, constraining them so that they satisfy $E_{k,y_{k+1}}\leq 0$, then the rest are free. This gives a count of $2^{n-\lceil\frac{n-1}{2}\rceil}\binom{n}{\lceil\frac{n-1}{2}\rceil}$ $y\in I_n$ with resepct to which 
\[
\int_{I_n} e^{t(d_H(x,y) -\E d_H(x,y))} \, d\mu(x) \leq e^{nt^2/2}.
\]
This completes the proof.
\end{Proof}



\subsubsection{The Small Variance Case}  \label{smallvariance}
In this section we further analyze the upper bound provided in Theorem \ref{errorbound}. We will observe that the form of the upper bound may be significantly simplified upon putting absolute values around the ``error-terms" $E_{k,x',y',t}$. As these may be negative, the bound we derive in this section will be worse than the bound in Theorem \ref{errorbound}. However, the bound is appealing due to its simplicity,  and it is well-suited for distributions which are approximately product distributions, i.e. ones whose conditional distributions (\ref{conditionalmeasure}) have a small variance. 

We begin with a lemma.

\begin{lemma}\label{productformulasmallvariance}
  For any sequence $b_j\in \R$, $j\geq 2$, we have 
 \begin{equation}\label{productformula}
 \prod_{j=2}^n (b_j+e^{t^2/2})
 =e^{(n-1)t^2/2} +e^{(n-2)t^2/2}b_2 +\sum_{k=2}^{n-1} e^{(n-k-1)t^2/2}b_{k+1}\prod_{j=2}^k (b_j+e^{t^2/2}).
 \end{equation}
\end{lemma}
 
 \begin{Proof}
 We proceed by induction. For the base case, $n=2$, the left and right hand side are both $b_2+e^{t^2/2}$.
 
 Now assume (\ref{productformula}) for some $n>2$. We have 
\begin{eqnarray*}
\prod_{j=2}^{n+1}(b_j+e^{t^2/2}) &=& (b_{n+1}+e^{t^2/2})\prod_{j=2}^n (b_j+e^{t^2/2})\\
&=& b_{n+1}\prod_{j=2}^n (b_j+e^{t^2/2})+e^{t^2/2}\prod_{j=2}^n (b_j+e^{t^2/2}).
\end{eqnarray*}
By the induction hypothesis,
$e^{t^2/2}\prod_{j=2}^n (b_j+e^{t^2/2})$ is 
\[ 
e^{nt^2/2} +e^{(n-1)t^2/2}b_2 +\sum_{k=2}^{n-1} e^{(n-k)t^2/2}b_{k+1}\prod_{j=2}^k (b_j+e^{t^2/2}).
\]
It follows that 
\[
\prod_{j=2}^{n+1}(b_j+e^{t^2/2})= e^{nt^2/2} +e^{(n-1)t^2/2}b_2 +\sum_{k=2}^n e^{(n-k)t^2/2}b_{k+1}\prod_{j=2}^k (b_j+e^{t^2/2}),
\]
as desired.
\end{Proof}

We apply this lemma to the error bound given in Theorem \ref{errorbound}.  
To simplify notation, for any $j\geq 2$ and $t>0$ we set 
\[
b_j := e^{-t\mu^{(j)}(x_j\neq y_j)}(e^t-1)\lVert \epsilon_{x',y_j}\rVert_{\infty}.
\]
Recall the definition of 
$E_{n,x',y',t}$, which we denote here by $E_n$, $n\geq 1$:
\[
E_n=E_{n,x',y',t} =
 e^{-t\mu^{(n+1)}(x_{n+1}\neq y_{n+1})}(1-e^t) \int_{I_{n}}  a_{n}(x',y',t)\epsilon_{x',y_{n+1}}\, d\mu_{n}(x').
  \] 

\begin{prop}\label{boundonerror}
 We have
 \[
 |E_n| \leq \begin{cases} 
 b_{n+1}e^{t^2/2} \prod_{j=2}^n(b_j+e^{t^2/2}),& n\geq 2 \\
 b_2e^{t^2/2}, & n=1
 \end{cases}
 \]
\end{prop}

\begin{Proof}
The $n=1$ case is immediate from Hoeffding's lemma, so assume the result for all $2\leq k\leq n-1$. We have 
\[
|E_n| \leq b_{n+1}\int_{I_n}
|a_n(x,y,t)| \, d\mu_n(x) 
\leq b_{n+1}\left(e^{nt^2/2}+\sum_{k=1}^{n-1} e^{(n-k-1)t^2/2}|E_k|\right)
\]
using Theorem \ref{errorbound} in the last inequality. Using the induction hypothesis, this is bounded above by
\[
 b_{n+1}e^{t^2/2}\left( e^{(n-1)t^2/2}+e^{(n-2)t^2/2}b_2 +\sum_{k=2}^{n-1} e^{(n-k-1)t^2/2}b_{k+1}\prod_{j=2}^k(b_j+e^{t^2/2})\right).
\]
Lastly, by Lemma \ref{productformulasmallvariance}, this is 
\[
b_{n+1}e^{t^2/2} \prod_{j=2}^n(b_j+e^{t^2/2}).
\]
This completes the proof.
\end{Proof}

In a similar fashion, we have
our ``small-variance" bound,
\begin{thm}\label{smallvariancebound}
For any $t>0$,
\[
\int_{I_n} e^{t(d_H(x,y) -\E_{\mu_n}d_H(x,y))} \, d\mu_n(x) \leq e^{t^2/2}\prod_{j=2}^n (b_j +e^{t^2/2}).
\]
\end{thm}
\begin{Proof}
We have 
\begin{eqnarray*}
\int_{I_n} e^{t(d_H(x,y) -\E_{\mu_n}d_H(x,y))} \, d\mu_n(x) &\leq& e^{nt^2/2} +\sum_{k=1}^{n-1} e^{(n-k-1)t^2/2}|E_k| \\
&\leq& e^{t^2/2} \prod_{j=2}^n(b_j+e^{t^2/2}),
\end{eqnarray*}
using Theorem \ref{errorbound} in the first inequality, and Lemma \ref{productformulasmallvariance}
and Proposition \ref{boundonerror} in the second (exactly as in the proof of Proposition \ref{boundonerror}).
\end{Proof}
  
\section{Extending the Small Variance Bound to the Class of Lipschitz Functions on $I_n$}\label{Lipschitzuniform}
In this section we consider generalizations of the ``small-variance" bound, Theorem \ref{smallvariancebound}, to the general class of Lipschitz functions. Again we denote by $\mu_n$ the joint distribution of the vector $(x_1,\dots,x_n)\in I_n$.
We assume a uniform bound on the conditional distributions of $\mu_n$, denoted by $\mu(x_k \mid x_1,\dots, x_{k-1})$,
namely
\begin{equation}\label{unifbound}
\mu(x_k \mid x_1,\dots, x_{k-1}) \leq c_k \leq 1
\end{equation}
with $c_k$ decreasing to $1/2$ as $k\to \infty$.
Assume that $\mu_1(0)=1/2$. 
\begin{thm}
Let $A\subset I_n$ be a non-empty set and let $\mu_n$ denote the joint distribution of the vector $(x_1,\dots, x_n)\in I_n$, satisfying (\ref{unifbound}) and $\mu_1(0)=1/2$. 
Then for all $t>0$, 
\begin{eqnarray*}
\int_{I_n} e^{td(x,A)} \, d\mu_n &\leq& \frac{1}{\mu_n(A)}\left(\frac{1}{2}+\frac{e^{t}+e^{-t}}{4}\right) \left(
\prod_{k=2}^nc_k^2\right)\big(2+e^t+e^{-t}\big)^{n-1} \\
&\leq& \frac{1}{\mu_n(A)}\left(4^{n-1}\prod_{k=2}^nc_k^2\right)e^{t^2n/4}
\end{eqnarray*}
\end{thm}

\begin{Proof} 
We proceed by induction, closely following the work of Talagrand \cite{talagrand}. When $n=1$, if $A$ consists of one point, then $\mu_1(A)=1/2$ and
\[
\int_{I_1}e^{td(x,A)} \, d\mu_1(x) = \frac{1}{2}+\frac{1}{2}e^t \leq 2\left(\frac{1}{2}+\frac{e^{t}+e^{-t}}{4}\right).
\]
If $A$ consists of two points, then $\mu_1(A)=1$ and 
\[
\int_{I_1}e^{td(x,A)} \, d\mu_1(x) = 1\leq \frac{1}{2}+\frac{e^{t}+e^{-t}}{4}.
\]

Suppose $n>1$ and assume the result for $n-1$. 
Define $I_n^1$ to be the set of all $x\in I_n$
with a one in the last coordinate, and similarly define $I_n^0$ to 
be the set of all $x\in I_n$ with a zero in the last coordinate.
Then of course $I_n^0$ and $I_n^1$ partition the space $I_n$.
Moreover, given $A\subset I_n$, we define 
\[
A_0= \{ x\in I_{n-1} : (x,0) \in A\}
\]
and
\[
A_1 = \{x\in I_{n-1} : (x,1) \in A\}.
\]
Note that $A = A_0\times \{0\} \cup A_1\times \{1\}$, and
\begin{eqnarray*}
\mu_n(A) &=& \sum_{x\in A_0} \mu_n(x,0) + \sum_{x\in A_1} \mu_n(x,1) \\
               &\leq&  \sum_{x\in A_0} c_n\mu_{n-1}(x) + \sum_{x\in A_1} c_n\mu_{n-1}(x) \\
               &=& c_n\mu_{n-1}(A_0)+c_n\mu_{n-1}(A_1).
               \end{eqnarray*}
Given $x\in I_n$, let $x'\in I_{n-1}$ denote $x$ with its last coordinate omitted. Then
\begin{equation}\label{distance_0}
d(x,A) =\min\{d(x', A_0), d(x',A_1) + 1\}, \;\; x\in I_n^0
\end{equation}
and
\begin{equation}\label{distance_1}
d(x,A) =\min\{d(x', A_1), d(x',A_0) + 1\}, \;\; x\in I_n^1.
\end{equation}

We now consider 
\[
\int_{I_n} e^{td(x,A)} \, d\mu_n  = \int_{I_n^0} e^{td(x,A)} \, d\mu_n +\int_{I_n^1} e^{td(x,A)} \, d\mu_n
\]
Define $C_n=\prod_{k=2}^nc_k^2$. Using (\ref{distance_0}) we have 
\begin{eqnarray*}
\int_{I_n^0} e^{td(x,A)} \, d\mu_n &=& 
\int_{I_n^0} \exp(\min\{td(x', A_0), td(x',A_1) + t\})\, d\mu_n  \\
&=& \int_{I_n^0} \min\{\exp(td(x', A_0)), e^t\exp(td(x',A_1))\}\, d\mu_n  \\
&=& \sum_{x'\in I_{n-1}}  \min\{\exp(td(x', A_0)), e^t\exp(td(x',A_1))\} \mu_n(x',0)\\
&\leq& c_n\sum_{x'\in I_{n-1}}  \min\{\exp(td(x', A_0)), e^t\exp(td(x',A_1))\} \mu_{n-1}(x')\\
&=& c_n\int_{I_{n-1}} \min\{\exp(td(x', A_0)), e^t\exp(td(x',A_1))\}\, d\mu_{n-1}(x') \\
&\leq& c_n \min\left\{\int_{I_{n-1}}\exp(td(x', A_0))\, d\mu(x'),
\int_{I_{n-1}}e^t\exp(td(x',A_1))\, d\mu(x')\right\} \\
&\leq& c_n\min\left\{ \frac{C_{n-1}c(t)^{n-1}}{\mu_{n-1}(A_0)}, e^t\frac{C_{n-1}c(t)^{n-1}}{\mu_{n-1}(A_1)}\right\}.
\end{eqnarray*}
Here the first inequality follows from (\ref{unifbound}), the second inequality follows
since the integral of the minimum does not exceed the minimum of the integrals,
and finally the third inequality follows from the induction hypothesis (and in the second to last line note that $\mu=\mu_{n-1}$). Note that 
we similarly derive
\[
\int_{I_n^1} e^{td(x,A)} \, d\mu_n  \leq 
c_n\min\left\{ \frac{C_{n-1}c(t)^{n-1}}{\mu_{n-1}(A_1)}, e^t\frac{C_{n-1}c(t)^{n-1}}{\mu_{n-1}(A_0)}\right\},
\]
and therefore
\[
\int_{I_n} e^{td(x,A)} \, d\mu_n  \leq
c_n\min\left\{ \frac{C_{n-1}c(t)^{n-1}}{\mu_{n-1}(A_0)}, e^t\frac{C_{n-1}c(t)^{n-1}}{\mu_{n-1}(A_1)}\right\}
\]
\[
+c_n\min\left\{ \frac{C_{n-1}c(t)^{n-1}}{\mu_{n-1}(A_1)}, e^t\frac{C_{n-1}c(t)^{n-1}}{\mu_{n-1}(A_0)}\right\}.
\]
Write the right hand side as
\[
\frac{C_{n-1}c(t)^{n-1}}{\mu_{n}(A)}\left(c_n\min\left\{ \frac{\mu_n(A)}{\mu_{n-1}(A_0)}, e^t\frac{\mu_n(A)}{\mu_{n-1}(A_1)}\right\}
+c_n\min\left\{ \frac{\mu_n(A)}{\mu_{n-1}(A_1)}, e^t\frac{\mu_n(A)}{\mu_{n-1}(A_0)}\right\} \right)
\]
Now denote
\[
a_0 = \frac{\mu_{n-1}(A_0)}{\mu_n(A)}, \;\; a_1=\frac{\mu_{n-1}(A_1)}{\mu_n(A)}.
\]
From the inequality $\mu_n(A) \leq c_n(\mu_{n-1}(A_0) + \mu_{n-1}(A_1))$, we 
see that $a_0+a_1\geq \frac{1}{c_n}$. 
To complete the proof we deduce the maximum
value of 
\[
c_n\min\{a_0^{-1},e^ta_1^{-1}\} +c_n\min\{a_1^{-1},e^ta_0^{-1}\}
\]
where $a_0+a_1 \geq \frac{1}{c_n}$. (These minimization problems
are not independent from one another.)

We first assume that $a_0+a_1 =\frac{1}{c_n}$. Set $c(t):= 2+e^t+e^{-t}$. We split the maximization
into cases:
\begin{itemize}
\item
If $a_0=0$, then $a_1=\frac{1}{c_n}$ and the value is $c_n^2+e^tc_n^2\leq c_n^2c(t)$
\item
Assume that $a_0,a_1>0$. If $a_1=a_0=\frac{1}{2c_n}$, the value is
$2c_n^2+2c_n^2=4c_n^2$.  We have $4c_n^2\leq c_n^2c(t)$, since
$2\leq e^t+e^{-t}$.
\item
Assume that $a_0,a_1>0$ and $a_0\neq a_1$. Without loss of generality
assume that $a_0>a_1$, so $a_0=\frac{1}{2c_n}+b$, $a_1=\frac{1}{2c_n}-b$,
for some $0<b<\frac{1}{2c_n}$. Then $a_0^{-1}<e^ta_1^{-1}$. If 
$a_1^{-1}<e^ta_0^{-1}$, the value is 
\[
c_n(1/2c_n+b)^{-1} + c_n(1/2c_n-b)^{-1} = \frac{1}{\frac{1}{4c_n^2}-b^2}.
\]
Increasing $b$ slightly we increase the value, a contradiction. 
A similar calculation reveals a contradiction if $a_1^{-1}>e^ta_0^{-1}$. Thus
at the maximum point we must have $a_1^{-1}=e^ta_0^{-1}$. This says
that $b=\frac{e^t-1}{2c_n(1+e^t)}$, and therefore
\[
a_0= \frac{e^t}{c_n(1+e^t)}
\]
and
\[
a_1= \frac{1}{c_n(1+e^t)}.
\]
This gives a value of 
\[
c_na_1^{-1}+c_na_0^{-1} =2c_n^2+c_n^2(e^t+e^{-t}) =c_n^2c(t).
\]
\end{itemize}
Finally, when $a_0+a_1\geq 1/c_n$, we must have $a_0+a_1=\frac{1}{b}$ for
some $b<c_n$, and running the above argument through again gives a maximum
value of $2b^2+b^2(e^t+e^{-t}) < 2c_n^2+c_n^2(e^t+e^{-t})$. Thus the maximum
is achieved when $a_0+a_1=1/c_n$, and the maximum value is $c_n^2c(t)$. 

It follows that,
\[
\int_{I_n} e^{td(x,A)} \, d\mu_n(x) \leq \frac{1}{\mu_n(A)}\left(\frac{1}{2}+\frac{e^{t}+e^{-t}}{4}\right) \left(
\prod_{k=2}^nc_k^2\right)\big(2+e^t+e^{-t}\big)^{n-1}.
\]
Lastly, noting that 
\[
\frac{1}{2}+\frac{e^{t}+e^{-t}}{4} \leq e^{t^2/4},
\]
which is clear upon power series expansion, we have 
\[
(2+e^t+e^{-t})^{n-1} \leq
4^{n-1} e^{t^2(n-1)/4},
\]
and therefore
\[
\int_{I_n} e^{td(x,A)} \, d\mu_n(x) \leq
\frac{1}{\mu_n(A)}\left(4^{n-1}\prod_{k=2}^nc_k^2\right)e^{t^2n/4}.
\]
The proof is complete.
\end{Proof}


\begin{thebibliography}{99} 

\bibitem{spencer}  N. Alon and J. Spencer, {\it The Probabilistic Method}, Wiley, (1990).

\bibitem{ball}  K. Ball, {\it An Elementary Introduction to Modern Convex Geometry}, Flavors of Geometry 31, (1997).

\bibitem{barvinok} A. Barvinok , {\it Math 710 : Measure Concentration}, (2005), {\tt https://dept.math.lsa.umich.}\\{\tt edu/~barvinok/total710.pdf}

\bibitem{gromov-milman}  M. Gromov and V.D. Milman, {\it A topological application of the isoperimetric inequality}, Amer. J. Math. 105, (1983), pp. 843-854.

\bibitem{Kontor}  L. Kontorovich and K. Ramanan, {\it Concentration inequalities for dependent random variables via the martingale method}, The Annals of Probability 36 (2008), pp. 2126-2158.

\bibitem{CM} M. Ledoux, {\it The Concentration of Measure Phenomenon}, American Mathematical Society (2005).

\bibitem{marton3}  K. Marton, {\it Measure concentration for Euclidean distance in the case of dependent random variables}, The Annals of Probability 32, (2004), pp. 2526-2544.

\bibitem{marton2} K. Marton, {\it Measure concentration and strong mixing}, Studia Scientiarum Mathematicarum Hungarica 40, (2003), pp. 95-113.

\bibitem{Marton} K. Marton, {\it A measure concentration inequality for contracting Markov chains}, Geometric and Functional Analysis 6, (1996), pp. 556-571.

\bibitem{asymp}  V.D. Milman and G. Schechtman, {\it Asymptotic Theory of Finite Dimensional Normed Spaces}, Lecture Notes in Mathematics. Springer Berlin Heidelberg, (1986).

\bibitem{milman} V.D. Milman, {\it A new proof of A. Dvoretzky's theorem on cross-sections of convex bodies}, Funkcional Anal i Prilozen (1971), pp. 28-37.

\bibitem{pemantle} R. Pemantle, {\it Towards a theory of negative dependence}, Journal of Mathematical Physics 41, (2000), pp. 1371-1390.

\bibitem{peres} R. Pemantle and Y. Peres, {\it Concentration of Lipschitz functionals of determinantal and other strong Rayleigh measures}, Combin. Probab. Comput. 23 (2014), pp. 140-160.

\bibitem{Samson}  P.M. Samson, {\it Concentration of measure inequalities for markov chains and $\Phi$-mixing processes}, The Annals of Probability 28, (2000), pp. 416-461. 

\bibitem{talagrand} M. Talagrand, {\it A new look at independence}, Annals of Probability 24, (1996), pp.1-34.

\bibitem{TT} T. Tao, {\it Topics in Random Matrix Theory}, Graduate Studies in Mathematics, American Mathematical Society, (2012).





\end{thebibliography}
\end{document}